\documentclass[12pt]{amsart}
\usepackage{amsmath,amssymb,amsfonts,amsthm}
\usepackage{mathrsfs}
\usepackage{hyperref}
\usepackage[margin=0.75in]{geometry}

\newtheorem*{claim}{Claim}
\newtheorem{thm}{Theorem}[section]
\newtheorem{prop}[thm]{Proposition}
\newtheorem{cor}[thm]{Corollary}
\newtheorem{lem}[thm]{Lemma}

\newtheorem{conj}[thm]{Conjecture}

\theoremstyle{definition}
\newtheorem{defn}[thm]{Definition}

\theoremstyle{remark}
\newtheorem{rmk}[thm]{Remark}
\newtheorem*{pf}{Proof}
\newtheorem{eg}[thm]{Example}

\newcommand\cO{{\mathcal{O}}}
\newcommand\bC{{\mathbb C}}
\newcommand\bN{{\mathbb N}}
\newcommand\bQ{{\mathbb Q}}
\newcommand\bR{{\mathbb R}}
\newcommand\bZ{{\mathbb Z}}
\newcommand\sF{{\mathscr{F}}}
\newcommand\sG{{\mathscr{G}}}
\newcommand\sH{{\mathscr{H}}}
\newcommand\rw{\longrightarrow}

\begin{document}

\title{ACC for foliated log canonical thresholds}


\author{Yen-An Chen}
\address{School of Mathematics, Korea Institute for Advanced Study, 85 Hoegi-ro, Dongdaemun-gu, Seoul 02455, Republic of Korea.}
\email{yachen@kias.re.kr}

\begin{abstract}
It is known that the set of log canonical thresholds (lcts) on any varieties with fixed dimension satisfies the ascending chain condition. Inspired by the foliated minimal model program, it is intriguing to study the foliated version of lcts and ask whether they have the similar property. We give an affirmative answer in the case of surfaces and threefolds. 
\end{abstract}

\maketitle

\section*{Introduction}
Singularities play an important role in studying higher dimensional birational geometry. 
In order to measure the singularities, we would like to find some invariants which behave well under birational maps. 
For example, the minimal log discrepancies (mlds for short) arise naturally in the study of birational geometry. 
In \cite{shokurov2004letters5}, Shokurov showed that two following conjectures on minimal log discrepancies implies the termination of flips, which is one of main open problems in the minimal model program. 

\begin{conj}[{\cite[Problem 5]{shokurov1988problems}}]\label{acc_conj}
    Fix a positive integer $n$ and a subset $I\subset [0,1]$ satisfying the descending chain condition. 
    Then the set 
    \[\operatorname{MLD}_n(I) := \{\operatorname{mld}(x,X,B)\vert\,(X\ni x,B) \textnormal{ is log canonical},\, \dim X = n,\, B\in I\}\]
    satisfies the ascending chain condition. 
    Here $\operatorname{mld}(x,X,B)$ is the minimal log discrepancy of the pair $(X,B)$ with centre $x$ and $B\in I$ means that the coefficients of $B$ belong to $I$. 
\end{conj}
\begin{conj}[{\cite[Conjecture 2.4]{ambro1999minimal}}]\label{lsc_conj}
    Let $X$ be a quasi-projective normal variety and $B$ be an $\bR$-Weil divisor such that $K_X+B$ is $\bR$-Cartier. 
    Then the function $a \colon |X| \to \bQ$, which sends a closed point $x\in X$ to the minimal log discrepancy of $(X,B)$ with centre $x$, is lower semi-continuous.
\end{conj}

Conjecture~\ref{acc_conj} is known to hold for surface pairs \cite{alexeev1993two}, toric varieties \cite{borisov1997minimal}, toric pairs \cite{ambro2006toric}, exceptional singularities \cite{han2019acc}, and recently terminal threefolds \cite{han2022acc}; while conjecture~\ref{lsc_conj} holds for $n\leq 3$ \cite{alexeev1993two, ambro1999minimal} and $X$ is toric \cite{ambro1999minimal}. 

Another interesting and important invariant is the log canonical threshold, which plays a vital role in the inductive approach for higher dimensional geometry. 
The set of log canonical thresholds also satisfies the ascending chain condition by the following theorem. 
\begin{thm}[{\cite[Theorem 1.1]{hacon2014acc}}]
    Fix $n\in\bN$, $I\subset [0,1]$, and a subset $J$ of positive real numbers. 
    If $I$ and $J$ both satisfy the descending chain condition, then the set
    \[\operatorname{LCT}_n(I,J) = \{\operatorname{lct}(X,\Delta;M)\vert\, (X,\Delta)\in\mathfrak{T}_n(I) \mbox{ and } M \in J\}\]
    satisfies the ascending chain condition where $M$ is $\bR$-Cartier and $\mathfrak{T}_n(I)$ is the set of log canonical pairs $(X,\Delta)$ where $X$ is a normal variety of dimension $n$ and $\Delta\in I$. 
\end{thm}
It was originally conjectured in all dimensions and proved for $n=2$ by Shokurov in \cite{shokurov1992three}. 
And then Alexeev proved it for $n=3$ in \cite{alexeev1993two}. 
Later, it was proved for smooth varieties in \cite{defernex2010acc} and for full generality in \cite{hacon2014acc}. 

Inspired by the work on the minimal model program on foliated surfaces and threefolds (see, for example, \cite{brunella2015birational, mcquillan2008canonical, spicer2020higher, cascini2021mmp, spicer2022local, cascini2025mmp} and the references therein), we are interested in the foliated version of the theorem above. 
Precisely, we define the foliated log canonical threshold as follows.
\begin{defn}\label{defn_fol_lct}
    Let $(\sF,\Delta)$ be a log canonical pair (see Definitions~\ref{defn_folpair} and \ref{defn_lc}) on a normal variety $X$.  
    \begin{enumerate}
        \item For an $\bR$-Cartier divisor $M\geq 0$, we define the \emph{log canonical threshold} of $M$ with respect to $(\sF,\Delta)$ to be 
        \[\operatorname{lct}(\sF,\Delta;M) := \sup\{t\in\bR\vert\,(\sF,\Delta+tM) \textnormal{ is log canonical}\}.\]
        \item Let $I\subset [0,1]$ and $J\subset\bR_{> 0}$. 
        We put $\operatorname{LCT}_{n,\,r}(I,J)$ to be the set of all real numbers $\operatorname{lct}(\sF,\Delta;M)$ where $(\sF,\Delta)$ is a log canonical pair such that $\sF$ is a foliation of rank $r$ on a variety $X$ of dimension $n$ and the coefficients of $\Delta$ (resp. $M$) belong to $I$ (resp. $J$). 
    \end{enumerate}
\end{defn}

In this paper, we show the following theorem.
\begin{thm}\label{main_acc}
    Fix $I\subset [0,1]$, a subset $J$ of positive real numbers, and two positive integers $r$ and $n$ with $r<n$. 
    If $I$ and $J$ both satisfy the descending chain condition, then the set $\operatorname{LCT}_{n,\,r}(I,J)$ satisfies the ascending chain condition provided that $n\leq 3$. 
\end{thm}

With some reformulations as in \cite{hacon2014acc}, we also show the following theorem. 

\begin{thm}\label{thm_coef_finite}
    Fix two positive integers $r$ and $n$ with $r<n\leq 3$ and a set $I\subset [0,1]$, which satisfies the descending chain condition. 
    Then there exists a finite subset $I_0\subset I$ with the following property:
    If $(X,\sF,D)$ is a triple such that 
    \begin{enumerate}
        \item $X$ is a normal variety of dimension $n$,
        \item $\sF$ is a foliation of rank $r$, 
        \item $(\sF,D)$ is log canonical,
        \item the coefficients of $D$ belong to $I$, and 
        \item there is a log canonical centre $Z$ which is contained in every component of $D$,
    \end{enumerate}
    then the coefficients of $D$ belong to $I_0$. 
\end{thm}

The proof has two steps. 
We first extract the log canonical centre and then apply adjunction to the foliated triple restricting onto some irreducible extracted divisor. 
The main difficulty is that, in general, there is little control on the singularities when applying adjunction. 
However, the model that extracts the log canonical centre has mild singularities, which gives us more control on singularities. 
Note that both steps rely on the foliated minimal model program and the resolution of foliation singularities, which are widely open when the dimension is at least $4$. 

The paper is organized as follows. 
In section 1, we fix some notations and recall some facts. 
In section 2, we show Theorem~\ref{thm_coef_finite} in the surface case. 
In the section 3 and 4, Theorem~\ref{thm_coef_finite} for threefolds is proved in the cases when the rank is 2 and 1, respectively. 
The last section is to show the set of log canonical thresholds satisfies the ascending chain condition (Theorem~\ref{main_acc}) from our finiteness theorem~\ref{thm_coef_finite}.

\section*{acknowledgements}
    The author would like to thank Jungkai Alfred Chen for his suggestions and encouragements. 
    He would also like to thank Paolo Cascini and Calum Spicer for reading the early draft. 
    He also wants to thank Jihao Liu for pointing out some unclarities in the proof of the first version and the anonymous referees for valuable comments.

\section{Preliminaries}
We will work over $\bC$.
\subsection{Foliations} 
\begin{defn}
    A \emph{foliation} $\sF$ on a normal variety $X$ is a coherent subsheaf $\sF$ of the tangent sheaf $T_X$ such that 
    \begin{enumerate}
        \item it is closed under the Lie bracket and 
        \item $\sF$ is saturated, that is, the quotient $T_X/\sF$ is torsion free.
    \end{enumerate}
\end{defn}

The \emph{rank} of $\sF$ is its rank as a sheaf. 
The \emph{co-rank} of $\sF$ is defined as $\dim X - \operatorname{rank}\sF$. 
We define the \emph{canonical divisor} $K_\sF$ for the foliation $\sF$ as $\cO_X(-K_\sF) = \det(\sF)$. 

Let $\sF$ be a foliation of rank $r$ on a normal variety $X$. 
We have a morphism $\Omega_X^{[1]}\to \sF^*$. 
Taking the double dual of $r$-th wedge product, we get a morphism 
$\varphi \colon \Omega_X^{[r]} \to \cO_X(K_\sF)$,
which yields a map $\phi \colon (\Omega_X^{[r]}\otimes\cO_X(-K_\sF))^{**} \to \cO_X$,
and we define the \emph{singular locus} of $\sF$ as the co-support of the image of $\phi$. 

\begin{defn}[Rational transform of foliations]
    Let $\sF$ be a foliation on a normal variety $X$. 
    \begin{enumerate}
        \item  Let $f \colon Y \dashrightarrow X$ be a dominant rational map and $U$ be an open subset of $X$ such that $f\vert_V \colon V \to U$ is an isomorphism where $V:= f^{-1}(U)$. 
        Note that $\sF\vert_U \subset T_U \cong T_V$. 
        By \cite[Exercise II.5.15]{hartshorne1977algebraic}, there is a coherent subsheaf $\sG$ of $T_Y$ such that $\sG\vert_V = \sF\vert_U \subset T_V $. 
        Then the \emph{pullback foliation} $f^*\sF$ is defined to be the saturation of $\sG$. 
        Note that $\sG$ is closed under the Lie bracket and, by \cite[Lemma 1.8]{hacon2021birational}, this definition is well-defined. 
        \item For a birational map $g \colon X \dashrightarrow Z$, we define the \emph{pushforward foliation} $g_*\sF$ as $(g^{-1})^*\sF$. 
    \end{enumerate}
\end{defn}

\begin{defn}
    Let $\sF$ be a foliation on a normal variety $X$ and a subvariety $W$ of $X$. 
    Let $U = X\setminus\big(\operatorname{Sing}(X)\cup\operatorname{Sing}(\sF)\cup\operatorname{Sing}(W)\big)$. 
    \begin{enumerate}
        \item $W$ is \emph{tangent} to $\sF$ if $T_W\vert_U\rw T_X\vert_U$ factors through $\sF\vert_U$. 
        Otherwise we say $W$ is \emph{transverse} to $\sF$. 
        \item $W$ is \emph{invariant} if $\sF\vert_U\rw T_X\vert_U$ factors through $T_W\vert_U$. 
    \end{enumerate}
\end{defn}

\begin{defn}[Foliated pair]\label{defn_folpair}
    A \emph{foliated pair} $(\sF,\Delta)$ on a normal variety $X$ consists of a foliation $\sF$ and an effective $\bR$-divisor $\Delta$ on $X$ such that $K_\sF+\Delta$ is $\bR$-Cartier. 
    
    The assumption that $K_\sF+\Delta$ is $\bR$-Cartier is not necessary when $X$ is a normal surface. 
    In this case, the pullback of any Weil divisor can be defined using Mumford's intersection theory on its minimal resolution. 
\end{defn}

\begin{defn}[{\cite[Definition I.1.5]{mcquillan2008canonical}}]\label{defn_lc}
    Let $(\sF,\Delta)$ be a foliated pair on a normal variety $X$ and $f \colon Y\to X$ be a proper birational morphism from a normal variety $Y$. 
    For any divisor $E$ on $Y$, we define the \emph{discrepancy} of $(\sF,\Delta)$ along $E$ to be $a(E,\sF,\Delta) = \operatorname{ord}_E(K_{f^*\sF}-f^*(K_\sF+\Delta))$. 
    
    We say $(X,\sF,\Delta)$ is \emph{terminal} (resp. \emph{canonical}) if $a(E,\sF,\Delta)> 0$ (resp. $\geq 0$) for every exceptional divisor $E$ over $X$ and 
    $(X,\sF,\Delta)$ is \emph{log terminal} (resp. \emph{log canonical}) if $a(E,\sF,\Delta)> -\varepsilon(E)$ (resp. $\geq -\varepsilon(E)$) for every divisor $E$ over $X$ where $\varepsilon(E)$ is defined to be $0$ if $E$ is $f^*\sF$-invariant, and $1$ otherwise. 
\end{defn}

\begin{defn}
    Let $(\sF,\Delta)$ be a foliated pair on a normal variety $X$. 
    A subvariety $Z$ of $X$ is called a \emph{log canonical centre} for $(\sF,\Delta)$ if $(\sF,\Delta)$ is log canonical at the generic point of $Z$ and there is a divisor $E$ over $X$ with discrepancy $a(E,\sF,\Delta) = -\varepsilon(E)$ and with centre $Z$ on $X$. 
\end{defn}

\begin{thm}[{\cite{miyaoka1987deformations}, \cite[Theorem 9.0.2]{SB1992miyaoka}, or \cite{bogomolov2016rational}}]\label{BB}
    Let $\sF$ be a foliation on a normal variety $X$ of dimension $n$ and $H_1$, $\ldots$, $H_{n-1}$, and $A$  be ample divisors. 
    Let $C$ be a general intersection of element $D_i\in |m_iH_i|$ for some $m_i\gg  0$. 
    Suppose $K_\sF$ is $\bR$-Cartier. 
    If $K_\sF\cdot C<0$, then through a general point of $C$ there is a rational curve $\Sigma$ tangent to $\sF$ such that 
    \[A\cdot\Sigma\leq 2n\frac{A\cdot C}{-K_\sF\cdot C}.\]
\end{thm}

\subsection{Indices on foliated surfaces}
Most definitions in this subsection follow from \cite{brunella2015birational} with some generalizations. 

Let $p\in\operatorname{Sing}(\sF)\setminus\operatorname{Sing}(X)$. That is, $p$ is a smooth point on $X$ but a singular point of the foliation. 
Let $v$ be the vector field around $p$ generating $\sF$. 
Since $p\in\operatorname{Sing}(\sF)$, we have $v(p)=0$. 
Then we can consider the eigenvalues $\lambda_1$, $\lambda_2$ of $(\textnormal{D}v)\vert_p$, which do not depend on the choice of $v$. 

\begin{defn}\label{semi_reduced}
    If one of the eigenvalues is non-zero, say $\lambda_2$, then we define the eigenvalue of the foliation $\sF$ at $p$ to be 
    \[\lambda := \frac{\lambda_1}{\lambda_2}.\] 
    For $\lambda\neq 0$, this definition is unique up to reciprocal $\lambda \sim \frac{1}{\lambda}$. 
    
    If $\lambda = 0$, then $p$ is called a \emph{saddle-node}; otherwise, we say $p$ is \emph{non-degenerate}. 
    If $\lambda\notin\bQ^+$, then $p$ is called a \emph{reduced} singularity of $\sF$. 
\end{defn}

Reduced singularities arise naturally. 
Indeed, blowing up a smooth foliation point will introduce a reduced singularity with $\lambda=-1$. 

\subsubsection{Non-invariant curves}
We first consider the non-invariant curves and define the tangency order for them. 
\begin{defn}
    Let $(X,\sF)$ be a normal foliated surface and $C$ be a \emph{non-invariant} reduced curve. 
    Let $p\in C\setminus\operatorname{Sing}(X)$ and $v$ be the vector field generating $\sF$ around $p$. 
    Let $f$ be the local defining function of $C$ at $p$.  
    We define the \emph{tangency order} of $\sF$ along $C$ at $p$ to be 
    \[\operatorname{tang}(\sF,C,p) := \dim_\bC\frac{\cO_{X,p}}{\langle f,v(f)\rangle}.\]
    Note that $\operatorname{tang}(\sF,C,p)\geq 0$ and is independent of the choices of $v$ and $f$. 
    Moreover, if $\sF$ is transverse to $C$ at $p$, then $\operatorname{tang}(\sF,C,p) = 0$. 
    Therefore, if $C$ is compact, then we can define 
    \[\operatorname{tang}(\sF,C) := \sum_{p\in C}\operatorname{tang}(\sF,C,p).\]
\end{defn}

\begin{thm}[{Adjunction for non-invariant divisors, see \cite{brunella1997feuilletages}, \cite[Proposition 2.2]{brunella2015birational}}]\label{adjunct_non_inv}
    Let $\sF$ be a foliation on a smooth surface $X$. 
    Let $C$ be a non-invariant irreducible compact curve on $X$ and $C^\nu$ be the normalization of $C$.  
    Then there is an effective divisor $\Delta$ on $C^\nu$ such that $(K_\sF+C)\vert_{C^\nu} = \Delta$ with $\deg\Delta = \operatorname{tang}(\sF,C)$.
\end{thm}

\begin{cor}\label{cor_adjunct_non_inv}
    Let $C$ be a non-invariant compact curve on a normal foliated surface $(X,\sF)$. 
    If $C$ is contained in the smooth locus of $X$, then we have $(K_\sF+C)\cdot C\geq 0$. 
\end{cor}

\subsubsection{Invariant curves}
Now we study the invariant curves. 
\begin{defn}
    Let $(X,\sF)$ be a normal foliated surface and $C$ be an invariant curve. 
    Let $p\in C\setminus\operatorname{Sing}(X)$ and $\omega$ be a $1$-form generating $\sF$ around $p$. 
    If $C$ is an invariant curve and $f$ is the local defining function of $C$ at $p$, then we can write 
    \[g\omega = h\textnormal{d}f+f\eta\] 
    where $g$, $h$ are holomorphic functions, $\eta$ is a holomorphic $1$-form, and $h$, $f$ are relatively prime functions. 
    
    We define the index $\operatorname{Z}(\sF,C,p)$ to be the vanishing order of $\frac{h}{g}\vert_C$ at $p$. 
    This definition is independent of the choices of $f$, $g$, $h$, $\omega$, and $\eta$. 
    (For a reference, see \cite[page 15 in Chapter 2 and page 27 in Chapter 3]{brunella2015birational}.)
\end{defn}

Note that if $p\notin\operatorname{Sing}(\sF)$, then $\operatorname{Z}(\sF,C,p)=0$. 
Therefore, if $C$ is compact, then we can define 
\[\operatorname{Z}(\sF,C) := \sum_{p\in C} \operatorname{Z}(\sF,C,p)\]
where the sums are taken over only finitely many points. 

\begin{thm}[{Adjunction for invariant divisors, see \cite{brunella1997feuilletages}, \cite[Proposition 2.3]{brunella2015birational}}]\label{adjunct_inv}
    Let $\sF$ be a foliation on a smooth projective surface $X$. 
    Let $C$ be an invariant irreducible curve on $X$ and $C^\nu$ be its normalization. 
    Then there is an effective divisor $\Delta$ on $C^\nu$ such that $K_\sF\vert_{C^\nu} = K_{C^\nu}+\Delta$ with $\deg\Delta = \operatorname{Z}(\sF,C) + \deg\operatorname{Diff}_C(0)$ where $\operatorname{Diff}_C(0)$ is the different with $(K_X+C)\vert_{C^\nu} = K_{C^\nu}+\operatorname{Diff}_C(0)$. 
    In particular, we have $K_\sF\cdot C = \operatorname{Z}(\sF,C) + 2p_a(C)-2$ where $p_a(C)$ is the arithmetic genus of $C$. 
\end{thm}

\begin{lem}\label{fib_neg}
    Let $(X,\sF)$ be a normal foliated surface. 
    Suppose there is a projective morphism $f \colon X \to  C$ where $C$ is a curve. 
    Let $F$ be the general fibre of $f$. 
    \begin{enumerate}
        \item If $K_\sF\cdot F<0$, then $F$ is a smooth invariant rational curve with $K_\sF\cdot F = -2$. 
        \item If $K_X\cdot F<0$, then $F$ is a smooth rational curve with $K_X\cdot F = -2$. 
    \end{enumerate}
\end{lem}
\begin{pf}
    Since $F$ is a general fibre, we may assume that $F$ contains no singularity of $X$ and $\sF$. 
    \begin{enumerate}
        \item If $F$ is not invariant, then by adjunction for non-invariant divisors (Corollary~\ref{cor_adjunct_non_inv}), we have 
        \[0>K_\sF\cdot F = (K_\sF+F)\cdot F \geq 0,\]
        which is impossible. 
        Thus, $F$ is invariant. 
        Now by adjunction for invariant divisors (Theorem~\ref{adjunct_inv}), we have 
        \[0>K_\sF\cdot F = \operatorname{Z}(\sF,F) +2p_a(F)-2 = 2p_a(F)-2.\]
        Therefore, $p_a(F)=0$ and $K_\sF\cdot F = -2$. 
        \item By adjunction, we have 
        \[0>K_X\cdot F = (K_X+F)\cdot F = 2p_a(F)-2.\]
        Therefore, $p_a(F)=0$ and $K_X\cdot F = -2$. 
    \end{enumerate}
    \qed
\end{pf}

\subsection{Ascending/Descending chain condition}
Let $I$ be a set of real numbers. 
A sequence $\{a_i\}_{i=1}^\infty$ in $I$ is \emph{increasing} (resp. \emph{strictly increasing} resp. \emph{decreasing} resp. \emph{strictly decreasing}) if $a_i\leq a_{i+1}$ (resp. $a_i<a_{i+1}$ resp. $a_i\geq a_{i+1}$ resp. $a_i> a_{i+1}$) for all $i\in\bN$. 

We define 
\[I_+ := \{0\}\cup\Bigg\{j\, :\, j=\sum_{p=1}^\ell i_p \textnormal{ for some } i_1, \ldots, i_\ell\in I\Bigg\}\]
and 
\[\operatorname{D}(I) := \bigg\{a\leq 1\, : \, a = \frac{m-1+f}{m}, m\in\bN, f\in I_+\bigg\}.\]

We say $I$ satisfies the \emph{ascending chain condition} (resp. \emph{descending chain condition}) if it does not contain any infinite strictly increasing (resp. strictly decreasing) sequence in $I$. 
Equivalently, for any infinite increasing (resp. decreasing) sequence $\{a_i\}_{i=1}^\infty$ in $I$, there is a positive number $\ell$ such that $a_i=a_\ell$ for all $i\geq \ell$.

Also, if $I\subset [0,1]$ satisfies the descending chain condition, then the minimum among $I\setminus\{0\}$ exists, which is denoted by $I_{\textnormal{min}}$. 

\begin{prop}\label{D(I)_dcc}
    Let $I\subset [0,1]$. 
    We have the following properties:
    \begin{enumerate}
        \item If $I$ satisfies the descending chain condition, then so do $I_+$ and $\operatorname{D}(I)$. 
        \item $\operatorname{D}(\operatorname{D}(I)) = \operatorname{D}(I)$. 
    \end{enumerate}
\end{prop}
\begin{pf}
    This is straightforward. For a reference, see \cite[Lemma 4.4]{mckernan2004threefold}. \qed
\end{pf}

\begin{lem}\label{lem_sum_finite}
    Let $I$ be a set satisfying the descending chain condition and $\alpha$ be a positive real number. 
    Then there is a finite subset $I_0\subset I$ depending only on $I$ and $\alpha$ such that if $\sum_{j=1}^\ell k_j i_j = \alpha$ for some $k_j\in\bN$, $i_j\in I\setminus\{0\}$, and $\ell\in\bN$, then all $i_j\in I_0$. 
\end{lem}
\begin{pf}
    Note that 
    \[\alpha = \sum_{j=1}^\ell k_j i_j \geq \sum_{j=1}^\ell k_j I_{\textnormal{min}}.\]
    Thus the sum $\sum_{j=1}^\ell k_j$ is bounded, and so is the length $\ell$. 
    Without loss of generality, we may fix $\ell$ and $k_1$, $k_2$, $\ldots$, $k_\ell$. 
    Now we write 
    \[\alpha - k_1i_1 = \sum_{j=2}^\ell k_ji_j\]
    where the right hand side belongs to $I_+$ which satisfies the descending chain condition by Proposition~\ref{D(I)_dcc} while the left hand side belongs to $-(I_+-\alpha)=\{-(m-\alpha)\mid m\in I_+\}$ which satisfies the ascending chain condition as $I_+$ satisfies the descending chain condition. 
    
    The only set which satisfies both the ascending chain condition and the descending chain condition is the finite set. 
    Thus we have only finitely many possible $i_1$. 
    Similarly, there are only finitely many possible $i_j$. 
    Hence, there is a finite subset $I_0\subset I$ such that all $i_j$ belong to $I_0$ provided that $\sum_{j=1}^\ell k_j i_j = \alpha$. 
    \qed
\end{pf}

\subsection{Foliated dlt modification}
Most definitions in this subsection follow from \cite{cascini2021mmp}. 

\begin{defn}
    We say $z_1$, $\ldots$, $z_\ell\in\bC^*$ satisfy the \emph{non-resonant condition} if, for any non-negative integers $a_1$, $\ldots$, $a_\ell$ with $\sum_i a_iz_i=0$, we have $a_i=0$ for all $i = 1$, $\ldots$, $\ell$. 
\end{defn}

\begin{defn}\label{defn_simple_corank_one}
    Let $\sF$ be a co-rank one foliation on a smooth variety $X$ of dimension $n$. 
    We say a point $p\in X$ is a \emph{simple singularity} for $\sF$ if, in formal coordinates $x_1$, $\ldots$, $x_n$ around $p$, conormal sheaf $N_\sF^* := (T_X/\sF)^*$ is generated by a $1$-form which is in one of the following two forms, for some $1\leq r\leq n$:
    \begin{enumerate}
        \item There are $\lambda_1$, $\ldots$, $\lambda_r\in\bC^*$, which satisfy the non-resonant condition, such that 
        \[\omega = x_1\cdots x_r\cdot\sum_{i=1}^r\lambda_i\frac{\textnormal{d}x_i}{x_i}.\]
        \item There is an integer $k\leq r$ such that 
        \[\omega = x_1\cdots x_r\cdot\Bigg(\sum_{i=1}^k p_i\frac{\textnormal{d}x_i}{x_i}+\varphi(x_1^{p_1}\cdots x_k^{p_k})\sum_{i=2}^r\lambda_i\frac{\textnormal{d}x_i}{x_i}\Bigg)\]
        where $p_1$, $\ldots$, $p_k$ are positive integers without a common factor, $\varphi(s)$ is a formal power series which is not a unit, and $\lambda_2$, $\ldots$, $\lambda_r\in\bC^*$ satisfy the non-resonant condition.
    \end{enumerate}
\end{defn}

\begin{defn}
    Let $(\sF,\Delta)$ be a foliated pair of co-rank one on a normal variety $X$. 
    We say $(\sF,\Delta)$ is \emph{foliated log smooth} if the following hold:
    \begin{enumerate}
        \item $(X,\Delta)$ is log smooth,
        \item $\sF$ has simple singularities, and
        \item if $S$ is the support of the non-$\sF$-invariant components of $\Delta$, $p$ is a closed point, and $\Sigma_1$, $\ldots$, $\Sigma_k$ are (possibly formal) $\sF$-invariant divisors passing through $p$, then $S\cup\bigcup_{i=1}^k\Sigma_i$ is a normal crossing divisor at $p$. 
    \end{enumerate}
\end{defn}

\begin{prop}\label{fol_log_sm_tang_zero}
    Let $\sF$ be a (co-)rank one foliation on a normal surface $X$ and $C$ be a non-invariant curve on $X$. 
    Suppose $(\sF,C)$ is foliated log smooth, then $(K_\sF+C)\cdot C=0$. 
\end{prop}
\begin{pf}
    Note that $X$ is smooth since $(\sF,C)$ is foliated log smooth. 
    By adjunction for non-invariant divisors (Theorem~\ref{adjunct_non_inv}), we have $(K_\sF+C)\cdot C = \operatorname{tang}(\sF,C)\geq 0$. 
    Thus it suffices to show that tangency order of $\operatorname{tang}(\sF,C,p)=0$ for any $p\in C$. 
    
    If $p$ is a singularity for $\sF$, then there are two formal invariant divisors $\Sigma_1$ and $\Sigma_2$ passing through $p$. 
    Hence $C\cup\Sigma_1\cup\Sigma_2$ is not a normal crossing divisor at $p$, which contradicts the last requirement that $(\sF,S)$ is foliated log smooth. 
    
    Therefore, $p$ is a smooth point for $\sF$. 
    Let $\Sigma$ be the formal invariant divisor passing through $p$. 
    Since $C\cup\Sigma$ is a normal crossing divisor at $p$, we have $\operatorname{tang}(\sF,C,p)=0$. 
    \qed
\end{pf}

\begin{defn}
    Let $(\sF,\Delta)$ be a foliated pair of co-rank one on a normal variety $X$. 
    A \emph{foliated log resolution} of $(\sF,\Delta)$ is a birational morphism $\pi \colon Y \to X$ such that 
    \begin{enumerate}
        \item $\operatorname{Exc}(\pi)$ is a divisor and 
        \item $(\sG,\pi_*^{-1}\Delta+E)$ is foliated log smooth where $\sG$ is the pullback foliation on $Y$ and $E$ is the sum over all $\pi$-exceptional divisors.
    \end{enumerate}
\end{defn}

\begin{rmk}
    When $X$ is a surface, the existence of foliated log resolution follows from Seidenberg's result \cite{seidenberg1968reduction}. 
    For the threefolds $X$, such a resolution exists by \cite{cano2004reduction}. 
\end{rmk}

\begin{defn}
    We say a foliated pair $(\sF,\Delta)$ of co-rank one on a normal variety $X$ is \emph{foliated dlt} if 
    \begin{enumerate}
        \item every irreducible component of $\Delta$ is generically transverse to $\sF$ and has coefficient at most one, and 
        \item there is a foliated log resolution $\pi \colon Y\to X$ of $(\sF,\Delta)$ which only extracts divisor $E$ of discrepancy $-\varepsilon(E)$.
    \end{enumerate}
\end{defn}

\begin{defn}
    Let $(\sF,\Delta)$ be a foliated pair of co-rank one on a normal projective variety $X$. 
    We call a birational projective morphism $\pi \colon Y \to X$ is a \emph{foliated dlt modification} if $(\sG,\pi_*^{-1}\Delta+\sum\varepsilon(E_i)E_i)$ is foliated dlt and 
    \[K_\sG+\pi_*^{-1}\Delta+\sum\varepsilon(E_i)E_i+F = \pi^*(K_\sF+\Delta)\]
    for some effective $\pi$-exceptional $\bQ$-divisors $F$ on $Y$ where $\sG$ is the pullback foliation on $Y$ and the sum is over all $\pi$-exceptional divisors.
\end{defn}

\begin{thm}[{\cite[Theorem 8.1]{cascini2021mmp}}]\label{F_dlt_modification}
    Let $(\sF,\Delta)$ be a foliated pair of rank two on a normal threefold.
    Then $(\sF,\Delta)$ admits a foliated dlt modification $\pi\colon Y \to X$ such that if $\sG$ is the pullback foliation on $Y$, then 
    \begin{enumerate}
        \item $Y$ is $\bQ$-factorial, 
        \item $Y$ has at worst klt singularities, and
        \item $(\sG,\Gamma:=\pi_*^{-1}D+E_{\textnormal{n-inv}})$ is foliated dlt with $K_\sG + \Gamma = \pi^*(K_\sF + D)$ where $E_{\textnormal{n-inv}}$ is the sum over all $\pi$-exceptional non-invariant divisors. 
    \end{enumerate}
    Moreover, if $(\sF,\Delta)$ is log canonical, then we may choose $\pi \colon Y \to X$ such that any log canonical centre for $(\sG,\Gamma)$ is contained in its codimension one log canonical centre. 
\end{thm}

\section{Foliations of rank one on surfaces}
In this section, we show Theorem~\ref{surface}, which is Theorem~\ref{thm_coef_finite} when $r=1$ and $n=2$. 

\begin{lem}\label{lc_center_cyclic}
    Let $(\sF,\Delta)$ be a foliated pair of rank 1 of a normal surface. 
    Suppose $(\sF,\Delta)$ is log canonical at a point $p$. 
    If $p\in\operatorname{Supp}(\Delta)$, then $\sF$ is terminal at $p$. 
\end{lem}
\begin{pf}
    Since $(\sF, \Delta)$ is log canonical at $p$, we also have $\sF$ is log canonical at $p$. 
    Let $E=\cup_iE_i$ be the exceptional divisor for the minimal resolution of $\sF$ at $p$. 
    Note that the (foliated) discrepancy $a(E_i,\sF,\Delta)$ is strictly less than $a(E_i,\sF)$ for all $i$ because $p\in\operatorname{Supp}(\Delta)$. 
    
    Suppose $\sF$ is not terminal at $p$, then, by \cite[Lemma 4.2]{chen2023log}, there is an irreducible exceptional divisor $E_i$ such that $a(E_i,\sF) = -\varepsilon(E_i)$. 
    But this implies that 
    \[a(E_i,\sF,\Delta) < a(E_i,\sF) = -\varepsilon(E_i)\]
    which contradicts that $(\sF,\Delta)$ is log canonical at $p$. 
    Therefore, $\sF$ is terminal at $p$. 
    \qed
\end{pf}

The following theorem gives more explicit properties for foliated dlt modifications on surfaces, whose existences are given in \cite[Theorem 8.1]{cascini2021mmp}. 

\begin{thm}\label{dlt_model}
    Let $(\sF,\Delta,p)$ be a germ of foliated pair on a normal surface $X$ with $p\in\operatorname{Supp}(\Delta)$. 
    Suppose $p$ is a log canonical centre of $(\sF,\Delta)$. 
    Then we have a foliated dlt modification $\pi \colon (Y,\sG,\Delta_Y) \to (X,\sF,\Delta)$ around $p$ such that the following hold:
    \begin{enumerate}
        \item The graph of the $\pi$-exceptional divisors $E = \cup_iE_i$ is a chain with $a(E_i,\sF,\Delta)=0$. 
        \item $K_\sG+\Delta_Y = \pi^*(K_\sF+\Delta)$ where $\Delta_Y = \pi_*^{-1}\Delta$ is the proper transform of $\Delta$. 
        \item There is at most one terminal foliation point $q$ on $E$. 
        Moreover $q$ is only on the irreducible component $E_1$ which is one end of the chain $E$. 
        \item $\Delta_Y$ only meets $E_1$. 
    \end{enumerate}
\end{thm}

\begin{pf}
    Let $\pi\colon (Y,\sG,\Delta_Y)\to (X,\sF,\Delta)$ be a foliated dlt modification around $p$. 
    We have $\pi^*(K_\sF+\Delta)=K_\sG+\Delta_Y+\sum\iota(E_i)E_i$ where $E_i$ are prime exceptional divisors. 
    Note that $p$ is terminal for $\sF$, so all exceptional divisors are invariant, and thus all $\iota(E_i)=0$. 
    Hence, we have (ii). 
    
    Let $E=\sum E_i$ and $q$ be any point on $\operatorname{Supp}(E)$. 
    If $(\sG,\Delta_Y)$ is foliated log smooth at $q$, then $(Y,E)$ is log smooth at $q$. 
    So there are at most two irreducible components of $E$ passing through $q$. 
    On the other hand, if $(\sG,\Delta_Y)$ is not foliated log smooth at $q$, then it is terminal at $q$, otherwise $q$ is an lc centre for $(\sG,\Delta_Y)$, which contradicts \cite[Lemma 3.8]{cascini2021mmp}. 
    In particular $\sG$ is terminal at $q$. 
    Thus, there is exactly one irreducible component of $E$ passing through $q$. 
    Hence, the dual graph of $E$ is a chain, which gives (i). 
    Thus, we re-label $E_i$ so that $E_i\cdot E_{i+1}=1$ for $i=1$, $\ldots$, $m-1$ where $m$ is the number of irreducible components of $E$. 
    
    Note that, on $E_i$ for $1<i<m$, there is no terminal point for $\sG$ or reduced singularity other than $E_i\cap E_{i-1}$ and $E_i\cap E_{i+1}$. 
    Indeed, if there is some terminal point for $\sG$ or other reduced singularity, then $0=\pi^*(K_\sF+\Delta)\cdot E_i = (K_\sG+\Delta_Y)\cdot E_i\geq K_\sG\cdot E_i \geq 1-\frac{1}{r}\geq \frac{1}{2}$ for some $r\geq 2$, which is impossible. 
    
    However, by \cite[Theorem 3.4]{brunella2015birational}, there is a reduced singularity on $E$ other than the intersection points of irreducible components of $E$, which will be on either $E_1$ or $E_m$. 
    Let us say it is on $E_m$. 
    Thus, as $K_\sG\cdot E_i=0$ for $i\geq 2$, we have $\Delta_Y\cdot E_i=0$ for all $i\geq 2$. 
    Therefore, $\Delta_Y$ only meets $E_1$. 
    This shows (iv). 
    
    Regarding (iii), the remaining is to show there is at most one terminal foliation point on $E_1$. 
    Note that $\Delta\neq 0$ as $p\in\operatorname{Supp}(\Delta)$. 
    Then $\Delta_Y\cdot E\neq 0$ and thus, $\Delta_Y\cdot E_1 > 0$. 
    Suppose there are more than two terminal points on $E_1$, then $0=(K_\sG+\Delta_Y)\cdot E_1>K_\sG\cdot E_1 \geq 1-\frac{1}{r_1}-\frac{1}{r_2}$ for some $r_1$, $r_2\geq 2$, which is impossible. 
    \qed
\end{pf}

\begin{thm}\label{surface}
    Let $I\subset [0,1]$ be a set satisfying the descending chain condition. 
    Then there exists a finite subset $I_0\subset I$ with the following property:
    If $(X,\sF,D)$ is a triple such that 
    \begin{enumerate}
        \item $X$ is a normal surface,
        \item $\sF$ is a foliation of rank $1$, 
        \item $(\sF,D)$ is log canonical,
        \item the coefficients of $D$ belong to $I$, and 
        \item there is a log canonical centre $Z$ which is contained in every component of $D$,
    \end{enumerate}
    then the coefficients of $D$ belong to $I_0$. 
\end{thm}
\begin{pf}
    If $Z$ is a divisor, then $D=Z$ with coefficient one. 
    In other words, if this case can happen, then $1\in I$ and thus, we require $1\in I_0$. 
    
    If $Z$ is a point, by Theorem~\ref{dlt_model}, we have a foliated dlt modification 
    $\pi \colon (Y,\sG,D_Y) \to (X,\sF,D)$ 
    with $K_\sG+D_Y = \pi^*(K_\sF+D)$ where $\sG$ is the pullback foliation and $D_Y$ is the proper transform of  $D$. 
    Let $E_1$ be an irreducible component of the exceptional divisor which meets $D_Y$. 
    Then we have 
    \[0 = \pi^*(K_\sF+D)\cdot E_1 = (K_\sG+D_Y)\cdot E_1 = \frac{-1}{m}+\sum_j\frac{k_jd_j}{m}\]
    where $m$ is the index of the possible singularity on $E_1$ and $k_j\in\bN$ for all $j$. 
    Here $k_j$ is the multiple of $m$ when the corresponding irreducible component does not pass through the only singularity on $E_1$. 
    Thus we have $\sum_jk_jd_j  = 1$. 
    Since all $d_j$ belong to $I$ which satisfying the descending chain condition, by Lemma~\ref{lem_sum_finite}, we have a finite subset $I_0\subset I$ such that $d_j\in I_0$ for all $j$. 
    \qed
\end{pf}

Using some ideas in \cite[Theorem 5.3]{alexeev1993two}, we also show the similar result for the numerically trivial pairs. 
\begin{thm}\label{global_acc_fol}
    Let $I\subset [0,1]$ be a set satisfying the descending chain condition. 
    Then there is a finite subset $I_0\subset I$ with the following property:
    If $(X,\sF,B)$ is a triple such that 
    \begin{enumerate}
        \item $X$ is a normal surface, 
        \item $\sF$ is a foliation of rank one on $X$, 
        \item $(\sF,B)$ is log canonical, 
        \item $K_{\sF}+B$ is numerically trivial, and 
        \item the coefficients of $B$ belong to $I$, 
    \end{enumerate}
    then all coefficients of $B$ belong to $I_0$. 
\end{thm}
\begin{pf}
    Let $(X, \sF, B)$ be a triple satisfying (1)-(5). 
    We write $B = \sum_j b_jB_j$ where $b_j\in I$ and $B_j$ are irreducible components. 
    
    Let $I'=I\cup\{1\}$. 
    Let $p\colon X'\to X$ be a foliated dlt modification (see Theorem~\ref{F_dlt_modification}) and $K_{\sF'}+B'=p^*(K_\sF+B)$ where $\sF'$ is the pullback foliation of $\sF$ on $X'$ and $B'=p_*^{-1}B+\sum\varepsilon(E_i)E_i$. 
    Note that $\sF'$ has either reduced singularities or terminal singularities by Theorem~\ref{dlt_model}. 
    
    Let $A$ be an ample divisor on $X'$. 
    Since $\sH$ is of rank one, we have $\mu_{A,\,\textnormal{min}}(\sF') = \mu_A(\sF') = -K_{\sF'}\cdot A = B'\cdot A>0$.
    By \cite[Theorem 1.1]{campana2019foliations} or \cite[Proposition 2.2]{ou2023generic}, $\sF'$ is algebraically integrable. 
    Moreover, as $\sF'$ has only reduced or terminal singularities, through which there are at most two invariant curves passing, the algebraic integrability of $\sF'$ gives a fibration $\pi \colon X' \to T$ onto a curve $T$ and the general fibres are the leaves of $\sF'$. 
    
    Let $F$ be a general fibre of the fibration. 
    Since $F$ is invariant, by adjunction for invariant divisors (Theorem~\ref{adjunct_inv}), we have $K_{\sF'}\cdot F = -2$ and thus 
    \[2 = -K_{\sF'}\cdot F = B'\cdot F = \sum_{j}b_j(B'_j\cdot F).\]
    Note that all $B'_j\cdot F\in\bN$ for all $j$ since all $B'_j$ are non-invariant and $F$ is invariant and general. 
    Therefore, by Lemma~\ref{lem_sum_finite}, there is a finite subset $I'_0\subset I'$ such that $b_j\in I'_0$ for all $j$. 
    Hence, $I_0=I'_0\cap I\subset I$ is a desired finite subset. 
    \qed
\end{pf}

\section{Foliations of rank two on threefolds}
In this section, we show Theorem~\ref{acc_fol_2}, which is Theorem~\ref{thm_coef_finite} when $r=2$ and $n=3$. 

\subsection{Adjunction}
The following Proposition slightly generalizes \cite[Lemma 3.18]{cascini2021mmp} by providing more information on the foliated differents. 

\begin{prop}\label{fol_adjunct}
    Let $\sF$ be a co-rank 1 foliation on a normal variety $X$ of dimension at most $3$, $S$ be a prime divisor, and $I$ be  a subset of $[0,1]$. 
    Suppose $(\sF,\varepsilon(S)S+\Delta)$ is foliated dlt where $\Delta$ is an effective $\bQ$-divisor with coefficients belonging to the set $I$. 
    Assume that $K_X$, $K_X+\Delta$, and $S$ are $\bQ$-Cartier. 
    
    Then there exists an effective $\bQ$-divisor $\Theta$ such that $(K_\sF+\varepsilon(S)S+\Delta)\vert_{S^\nu} = K_\sG+\Theta$ where $S^\nu$ is the normalization of $S$, $\sG$ is the restricted foliation to $S^\nu$, and moreover, the following statements hold:
    \begin{enumerate}
        \item If $\varepsilon(S)=1$, then $(\sG,\Theta)$ is log canonical and the coefficients of $\Theta$ belong to $\operatorname{D}(I)$.  
        \item If $\varepsilon(S)=0$, then $(S^\nu, \Theta' = \lfloor\Theta\rfloor_{\textnormal{red}}+\{\Theta\})$ is log canonical and the coefficients of $\Theta$ belong to $\operatorname{D}(I)\cup\bN$. 
        Here $\{\Theta\} = \Theta-\lfloor\Theta\rfloor$. 
        
        Moreover, any component $C$ of $\Theta - \Theta'$ is a log canonical centre for $(\sF,\Delta)$ with $(K_\sF+\Delta)\vert_{C^\nu} = K_{C^\nu} + \Xi$ where the coefficients of $\Xi$ belong to $\operatorname{D}(I)\cup\bN$. 
    \end{enumerate}
\end{prop}

\begin{pf}
    The existence of $\Theta$ follows from \cite[Lemma 3.18]{cascini2021mmp}. 
    \begin{enumerate}
        \item When $\varepsilon(S)=1$, we have $(\sG,\Theta)$ is foliated dlt by \cite[Lemma 3.18]{cascini2021mmp}. 
        In particular, it is log canonical. 
        Let $C$ be any prime divisor on $S$ with $m_C\Theta\neq 0$. 
        Note that $C$ is not $\sG$-invariant. 
        \begin{enumerate}
            \item If $C$ is a log canonical centre for $(\sF,S+\Delta)$, then $(\sF,S+\Delta)$ is foliated log smooth at the generic point of $C$. 
            Thus, by \cite[Proposition 3.13(4)]{cascini2023foliation}, the coefficient of $\operatorname{Diff}_S(\sF)$ along $C$ is $0$ or $1$ and therefore, the coefficient of $\Theta$ along $C$ is in $\operatorname{D}(I)$.  
            \item On the other hand, if $C$ is not a log canonical centre for $(\sF,S+\Delta)$, then it is terminal at the generic point of $C$. 
            Thus, the index one cover associated to $K_\sF$ is smooth in a neighborhood of generic point of $C$. 
            Hence, by \cite[Proposition 3.13(4)]{cascini2023foliation}, the coefficient of $\Theta$ along $C$ is $\frac{a+m-1+f}{m}$ 
            where $a$ is a non-negative integer, $m$ is the Cartier index of $K_\sF$ at the generic point of $C$, and $f\in I_+$. 
            As $(\sG,\Theta)$ is log canonical, we have the coefficient of $\Theta$ is in $\operatorname{D}(I)$. 
        \end{enumerate}
        
        \item When $\varepsilon(S)=0$, we have $(S^\nu,\Theta')$ is log canonical by \cite[Lemma 3.18]{cascini2021mmp}. 
        We recall, in the proof of \cite[Lemma 3.18]{cascini2021mmp}, that $\Theta'$ is the different of $(\widehat{X},S+T+\Delta)$ with respect to $S$ where $\widehat{X}$ is the formal completion of $X$ along $S$ and $T$ is the sum of all formal invariant divisors meeting $S$. 
        By \cite[Corollary 16.7]{kollar1992flips}, we have $\Theta'\in\operatorname{D}(I)$ and therefore, the coefficients of $\Theta$ belong to $\operatorname{D}(I)\cup\bN$. 
        
        Besides, each irreducible component $C$ of $\Theta-\Theta'$ is a log canonical centre for $(\sF,\Delta)$ as shown in the proof of \cite[Corolloary 3.20]{cascini2021mmp}. 
        
        By \cite[Lemma 3.22]{cascini2021mmp}, we have 
        \begin{equation}\label{eq_subadj}
        (K_\sF+\Delta)\vert_{C^\nu} = K_{C^\nu} + \Xi
        \end{equation}
        for some effective divisor $\Xi$ where $\nu\colon C^\nu\to C$ is the normalization of $C$. 
        We recall the derivation of the equation~(\ref{eq_subadj}) in \cite[Lemma 3.22]{cascini2021mmp}. 
        Let $S'$ be the strong separatrix at a general point of $C$. 
        Then we take the adjunction 
        \[(K_\sF+\Delta)\vert_{S'^\nu} = K_{S'^\nu}+C'+\Delta_{S'^\nu}\]
        where $S'^\nu$ is the normalization of $S'$ and $\Delta_{S'^\nu}\geq 0$ is a $\bQ$-divisor whose coefficients belong to $\operatorname{D}(I)\cup\bN$. 
        Note that $(S'^\nu,C'+\lfloor\Delta_{S'^\nu}\rfloor_{\textnormal{red}}+\{\Delta_{S'^\nu}\})$ is log canonical and $(K_{S'^\nu}+C'+\Delta_{S'^\nu})\vert_{C'^\nu} = K_{C^\nu}+\Xi$.
        
        From \cite[Lemma 3.22]{cascini2021mmp}, $\nu(P)$ is a log canonical centre of $(\sF,\Delta)$ for any point $P$ contained in the support of $\lfloor \Xi\rfloor$. 
        Then by \cite[Lemma 3.8]{cascini2021mmp}, $(\sF,\Delta)$ is foliated log smooth at $\nu(P)$. 
        Therefore, the coefficients of $\Xi$ belong to $\operatorname{D}(\operatorname{D}(I))\cup\bN = \operatorname{D}(I)\cup\bN$.
    \end{enumerate}
    \qed
\end{pf}

\subsection{ACC for LCT when \texorpdfstring{$n=3$}{n=3} and \texorpdfstring{$r=2$}{r=2}}
\begin{thm}\label{acc_fol_2}
    Fix a set $I\subset [0,1]$, which satisfies the descending chain condition. 
    Then there exists a finite subset $I_0\subset I$, depending only on $I$, with the following property:
    Suppose that  
    \begin{enumerate}
        \item $X$ is a normal threefold,
        \item $\sF$ is a foliation of rank $2$ on $X$, 
        \item $(\sF,D)$ is log canonical,
        \item the coefficients of $D$ belong to $I$, and 
        \item there is a log canonical centre $Z$ which is contained in every component of $D$. 
    \end{enumerate}
    Then the coefficients of $D$ belong to $I_0$. 
\end{thm}
\begin{pf}
We may assume that $Z$ is maximal with respect to the inclusion. 
If $Z$ is a divisor, then $Z = D$ with coefficient one. 
Thus we require $1\in I_0$. 

Now suppose $Z$ is not a divisor, we have, by Theorem~\ref{F_dlt_modification}, a foliated dlt modification $\pi \colon Y \to X$ such that, if $\sG$ is the pullback foliation on $Y$, then we have the following properties:
\begin{enumerate}
\item $Y$ is $\bQ$-factorial. 
\item $Y$ has at worst klt singularities. 
\item $(\sG,\Gamma:=\pi_*^{-1}D+\sum_{i}\varepsilon(E_i)E_i)$ is foliated dlt with $K_{\sG} + \Gamma = \pi^*(K_{\sF} + D)$ where the sum is over all $\pi$-exceptional divisors. 
\item Any log canonical centre for $(\sG,\Gamma)$ is contained in its codimension one log canonical centre. 
\end{enumerate}

Let $D_1$ be one of irreducible components of $D$ and $d_1$ be the coefficient of $D_1$ in $D$. 
We will show that $d_1$ belongs to a finite subset $I_0\subset I$, depending only on $I$. 

Note that the proper transform $\pi_*^{-1}D_1$ has coefficient $d_1$ in $\Gamma$. 
Without loss of generality, we may assume that $E_1$ intersects $\pi_*^{-1}D_1$. 
By adjunction for foliated threefolds (Proposition~\ref{fol_adjunct}), we have 
\[(K_\sG + \Gamma)\vert_{E_1^\nu} = K_{\sH}+\Theta\]
and $\Theta\in\operatorname{D}(I)$ where $\sH$ is the restricted foliation on $E_1^\nu$. 
Since $\pi_*^{-1}D_1$ meets $E_1$, there is an irreducible component $\Theta_1$ of $\Theta$ whose coefficient has the form 
\begin{equation}\label{b_1}
b_1 := \frac{m-1+f+kd_1}{m}
\end{equation}
where $m$, $k\in\bN$ and $f\in I_+\cap [0,1]\subset\operatorname{D}(I)$. 
Note also that $\Theta_1$ dominates $Z$ by the assumption that $Z$ is contained in $D_1$. 

We have four cases based on the invariance of $E_1$ and the dimension of $Z$: 

\subsubsection{\texorpdfstring{$E_1$}{E_1} is non-invariant} 
By adjunction (Proposition~\ref{fol_adjunct}), we have $(\sH,\Theta)$ is log canonical and the coefficients of $\Theta$ belong to $\operatorname{D}(I)$. 
\begin{enumerate}
\item If $Z$, which is the centre of $E_1$, is a point, then $K_{\sH}+\Theta$ is numerically trivial. 
As $\operatorname{D}(I)$ satisfies the descending chain condition, by Theorem~\ref{global_acc_fol}, there is a finite subset $J_1\subset \operatorname{D}(I)$, depending only on $I$, such that all coefficients of $\Theta$, in particular $b_1$, belong to $J_1$. 
By \cite[Lemma 5.2]{hacon2014acc}, there is a finite subset $I_1\subset I$, depending only on $I$, such that $d_1\in I_1$.

\item If $Z$ is a curve, then $\psi := \pi\vert_{E_1^\nu}\colon E_1^\nu \to Z$ is a fibration. 
Since $\Theta_1$ dominates $Z$, we have $\Theta_1\cdot F > 0$ for a general fibre $F$ of $\psi$. 
Therefore, $K_\sH\cdot F = -\Theta\cdot F < 0$ since $(K_\sH+\Theta)\cdot F = 0$. 
By Lemma~\ref{fib_neg}, $F$ is invariant with $K_\sH\cdot F = -2$. 
Now we write $\Theta = \sum_jb_j\Theta_j$ where $\Theta_j$ are distinct irreducible components, then 
\[2 = -K_\sH\cdot F = \Theta\cdot F = \sum_j(\Theta_j\cdot F)b_j.\]
Recall that $b_j\in\operatorname{D}(I)$ and $\Theta_1\cdot F >0$. 
Since $(\sH,\Theta)$ is log canonical, each $\Theta_j$ is non-invariant. 
Then $\Theta_j\cdot F$ is a positive integer for all $j$ because $F$ is general and invariant. 
By Lemma~\ref{lem_sum_finite} for the set $\operatorname{D}(I)$, which satisfies the descending chain condition, and $\alpha=2$, there is a finite subset $J_2\subset \operatorname{D}(I)$, depending only on $I$, such that $b_j\in J_2$ for all $j$. 
By \cite[Lemma 5.2]{hacon2014acc}, there is a finite subset $I_2\subset I$, depending only on $I$, such that $d_1\in I_2$.
\end{enumerate}

\subsubsection{\texorpdfstring{$E_1$}{E_1} is invariant} 
By adjunction (Proposition~\ref{fol_adjunct}), we have 
\[(E_1^\nu, \Theta' = \lfloor\Theta\rfloor_{\textnormal{red}}+\{\Theta\})\] 
is log canonical and the coefficients of $\Theta$ belong to $\operatorname{D}(I)\cup\bN$. 

\begin{enumerate}
\item If $Z$ is a curve, then $\psi := \pi\vert_{E_1^\nu}\colon E_1^\nu \to Z$ is a fibration. 
Since $\Theta_1$ dominates $Z$, we have $\Theta_1\cdot F > 0$ for a general fibre $F$ of $\psi$. 
So $K_{E_1^\nu}\cdot F = -\Theta\cdot F < 0$ because $(K_{E_1^\nu} + \Theta)\cdot F = 0$. 
By Lemma~\ref{fib_neg}, we have $K_{E_1^\nu}\cdot F = -2$.  
If we write $\Theta = \sum_jb_j\Theta_j$, then 
\[2 = \Theta\cdot F = \sum_j(\Theta_j\cdot F)b_j.\]
If there is a $j_0$ such that $b_{j_0}\geq 2$ and $\Theta_{j_0}\cdot F \neq 0$, then 
\[2 = \sum_j(\Theta_j\cdot F)b_j \geq b_1 + b_{j_0}>0+2=2,\]
which is impossible. 
Hence, we may assume that all $b_j\in\operatorname{D}(I)$. 
Therefore, $b_1\in J_2$ by the choice of $J_2$ and thus $d_1\in I_2$ by the choice of $I_2$. 

\item If $Z$ is a point, then $K_{E_1^\nu}+\Theta$ is numerically trivial. 
When $\Theta = \Theta'$, we have a log canonical pair $(E_1^\nu,\Theta)$ such that $K_{E_1^\nu}+\Theta$ is numerically trivial and the coefficients of $\Theta$ belong to $\operatorname{D}(I)$. 
By \cite[Theorem 1.5]{hacon2014acc} for the set $\operatorname{D}(I)$, there is a finite subset $J_3\subset \operatorname{D}(I)$ such that all coefficients of $\Theta$ belong to $J_3$. 
And thus, by \cite[Lemma 5.2]{hacon2014acc}, $d_1$ belongs to a subset $I_3\subset I$, which depends only on $I$. 

Recall that the coefficient of $\Theta_1$ is $b_1=\frac{m-1+f+kd_1}{m}\in\operatorname{D}(I)\cup \bN$. 
\begin{claim}
    $b_1\leq 1$. 
\end{claim}
\begin{pf}
    If not, then $\Theta_1$ is a component of $\Theta-\Theta'$, which is an lc centre of $(\sG,\Gamma)$ by Proposition~\ref{fol_adjunct}. 
    Then $(\sG,\Gamma)$ is foliated log smooth at the generic point of $\Theta_1$ by \cite[Lemma 3.8]{cascini2021mmp}. 
    As $\Theta_1$ comes from the intersection of $E_1$, which is invariant, and $\pi_*^{-1}D_1$, which is not invariant. 
    By foliated log smoothness, $\sG$ is smooth at the generic point of $\Theta_1$, and thus, the coefficients of $\Theta_1$ in $\Theta$ and $\Theta'$ are the same, which contradicts that $\Theta_1$ is a component of $\Theta-\Theta'$. 
    This completes the proof of the Claim. 
\end{pf}

Now suppose $\Theta\neq\Theta'$. 
The proof of this case will be divided into several steps. 
Write $\lfloor\Theta\rfloor = \sum n_iC_i$ and $\{\Theta\} = \sum b_i\Theta_i$ where $n_i\in\bN$, $n_1\geq 2$, and $b_i\in\operatorname{D}(I)$. 
\begin{enumerate}
    \item We first give a summary of the following steps. 
    we first show, in Step (b), the Theorem holds when there is a $C_i$ with $n_i\geq 2$ intersecting $\Theta_1$. 
    Thus, we may assume we are in the case that no $C_i$ with $n_i\geq 2$ intersects $\Theta_1$. 
    Then, in Step (c), we show the Cartier index of $K_{E_1^\nu}$ at any point on $C_i$ with $n_i\geq 2$ is bounded above by a constant depending only on $I$. 
    For the remaining steps, we will run $(K_{E_1^\nu}+\Theta'-b_1\Theta_1)$-minimal model program (MMP) and do a case-by-case study. 
    
    \item In this step, we show the Theorem holds when there is a curve $C_i$ with $n_i\geq 2$ meeting $\Theta_1$. 

    If there is a curve $C_i$ with $n_i\geq 2$ intersecting $\Theta_1$, then $C_i$ is a component of $\Theta-\Theta'$. 
    By Proposition~\ref{fol_adjunct}, we have $(K_\sG+\Gamma)\vert_{C^\nu} = K_{C^\nu}+\Xi$ for some divisor $\Xi$ and $\Xi\in\operatorname{D}(I)\cup\bN$. 

    Thus, 
    \[0 = (K_{E_1^\nu}+\Theta)\cdot C_i = (K_\sG+\Gamma)\cdot C_i = \deg(K_{C_i^\nu} + \Xi).\]
    Since $\Theta\cdot C_i\neq 0$, we have $\deg(K_{C_i^\nu}) = -\deg\Xi < 0$. 
    Hence $\deg\Xi = -\deg(K_{C_i^\nu})=2$. 
    Therefore, if we write $\Xi = \sum_s u_sp_s$ where $p_s$ are distinct points on $C_i^\nu$, then we have $\sum_s u_s = 2$. 
    Recall $\Theta_1\cdot C_i\neq 0$ and $u_s\in\operatorname{D}(I)\cup\bN$. 
    If there is one $u_s$ larger than 1, then $\Xi=2p$. 
    By \cite[Lemma 3.22]{cascini2021mmp}, $p$ is an lc centre of $(\sG,\Gamma)$. 
    So $(\sG,\Gamma)$ is foliated log smooth at $p$, and thus the coefficient of $p$ in $\Xi$ is at most 1 by the construction of $\Xi$, which gives a contradiction. 

    Therefore, all $u_s$ belong to $\operatorname{D}(I)$. 
    From assumption that $C_i$ intersects $\Theta_1$, one of $u_s$, say $u_1$, has the form 
    \[u_1 = \frac{n-1+\ell b_1+g}{n}\in (0,1]\] 
    where $n$, $\ell \in\bN$ and $g\in I_+\cap [0,1]\subset\operatorname{D}(I)$. 
    Replacing $b_1$ by the equation~(\ref{b_1}), we have 
    \[u_1 = \frac{r-1+\ell kd_1 + h}{r}\]
    for some $r\in \bN$ and $h\in\operatorname{D}(I)$. 
    By the choice of $J_2$ and $\sum_su_s=2$, we have $u_1\in J_2$ and thus $d_1\in I_2$ by the choice of $I_2$. 

    \item In this step, we show the Cartier index of $K_{E_1^\nu}$ at any point on $C_i$ with $n_i\geq 2$ is bounded above by a constant $M$ depending only on $I$. 

    To show this, we use $C_1$ to denote any curve $C_i$ with $n_i\geq 2$ to simplify the notation in the proof. 
    By Proposition~\ref{fol_adjunct}, $C_1$ is a log canonical centre for the foliated dlt pair $(\sG,\Gamma)$. 
    Then $(\sG,\Gamma)$ is foliated log smooth at the generic point of $C_1$ by \cite[Lemma 3.8]{cascini2021mmp}. 
    Thus, there are two separatrices passing through $C_1$: one is $E_1$, which is the weak separatrix by \cite[Corollary 3.20]{cascini2021mmp}, and the other is denoted as $S$, which is the strong separatrix.  

    Now, by adjunction, we put $(K_\sG+\Gamma)\vert_{S^\nu} = K_{S^\nu}+\Delta$. 
    By \cite[Corollary 3.20]{cascini2021mmp}, we have 
    \begin{equation}\label{adjunct_strong}
    \Big(K_{\widehat{Y}}+S+E_1+\sum S_i+\Gamma\Big)\vert_{S^\nu} = K_{S^\nu}+\Delta'
    \end{equation}
    where $\widehat{Y}$ is the formal completion of $Y$ along $S$, $\Delta' = \lfloor\Delta\rfloor_{\textnormal{red}}+\{\Delta\}$, and $S$, $E_1$, $S_i$ are all invariant divisors on $\widehat{Y}$. 
    Note that if $(\Delta-\Delta')\cdot C_1\neq 0$, then any point $p$ in the support of $(\Delta-\Delta')\cap C_1$ is a log canonical centre for $(S^\nu,\Delta')$ and thus a log canonical centre for $(\widehat{Y}, S+E_1+\sum S_i+\Gamma)$ by inversion of adjunction. 
    Thus, by \cite[Lemma 8.14]{spicer2020higher}, $p$ is a log canonical centre for $(\sG,\Gamma)$ and hence $(\sG,\Gamma)$ is foliated log smooth at $p$ by \cite[Lemma 3.8]{cascini2021mmp}. 
    Therefore, $(\Delta-\Delta')\cdot C_1\in\bN$. 
    Also, we have 
    \begin{align*}
    0 &= (K_{S^\nu}+\Delta)\cdot C_1 \\
    &\geq (K_{S^\nu}+C_1)\cdot C_1 + (\Delta-C_1-\{\Delta\})\cdot C_1\\ 
    &\geq -2 +  (\lfloor\Delta\rfloor - C_1)\cdot C_1
    \end{align*}
    where the equality holds since $C_1$ is contracted to a point $Z$ by assumption. 
    Note that the coefficient of $C$ in $\Delta$ is at least 2 for any irreducible component $C$ of the support of $\Delta-\Delta'$. 
    Then there is at most one irreducible component of the support of $\Delta-\Delta'$ intersecting $C_1$. 
    Besides, if such irreducible component exists, then it is unique and its coefficient in $\Delta$ is $2$. 
    So $(\Delta-\Delta')\cdot C_1\in\{0,1\}$. 

    Next we also notice that
    \begin{align*}
    -(\Delta-\Delta')\cdot C_1 &= (K_{S^\nu}+\Delta')\cdot C_1 - (K_{S^\nu}+\Delta)\cdot C_1 \\
    &= \Big(K_{\widehat{Y}}+S+E_1+\sum S_i+\Gamma\Big)\cdot C_1 - 0 \\
    &= (K_{E_1^\nu}+\Theta')\cdot C_1 \\
    &= (K_{E_1^\nu}+\Theta)\cdot C_1 - (\Theta-\Theta')\cdot C_1 \\
    &= - (\Theta-\Theta')\cdot C_1 
    \end{align*}
    where the second equality comes from the equation~(\ref{adjunct_strong}) and that $C_1$ is contracted to a point $Z$ and the last equality comes from $K_{E_1^\nu}+\Theta\equiv 0$.
    Therefore, we have 
    \[(\Theta-\Theta')\cdot C_1 = (\Delta-\Delta')\cdot C_1 = 0 \textnormal{ or } 1\] 
    and thus $(K_{E_1^\nu}+\Theta')\cdot C_1 = 0$ or $-1$. 
    By adjunction, we have 
    \[(K_{E_1^\nu}+\Theta')\vert_{C_1^\nu} = K_{C_1^\nu}+\Psi\]
    where $\Psi = \sum_s u_sp_s$ with $u_s\in\operatorname{D}(I)$. 
    More precisely, $u_s$ has the form $\frac{m-1+f}{m}$ where $m$ is the Cartier index of $K_{E_1^\nu}$ at $p_s$ and $f\in I_+$. 
    Note that if $u_s = 1$, then $p_s$ is a log canonical centre for $(C_1^\nu,\Psi)$. 
    So by inversion of adjunction and \cite[Lemma 8.14]{spicer2020higher}, $p_s$ is a log canonical centre for $(\sG,\Gamma)$. 
    Since $(\sG,\Gamma)$ is foliated dlt, it is foliated log smooth at $p_s$ by \cite[Lemma 3.8]{cascini2021mmp}. 
    Thus the Cartier index of $K_{E_1^\nu}$ at $p_s$ is $1$ if $u_s=1$. 

    Now since $\deg(K_{C_1^\nu}) + \sum_su_s = 0$ or $-1$, we have $\sum_su_s= 2$ or $1$.
    By Lemma~\ref{lem_sum_finite}, there is a finite subset $J_4\subset \operatorname{D}(I)$ such that each $u_s$ belongs to $J_4$. 
    We have seen above that the Cartier index of $K_{E_1^\nu}$ at $p_s$ is $1$ if $u_s=1$. 
    As $J_4$ is finite, the Cartier index of $K_{E_1^\nu}$ at $p_s$ is bounded above if $u_s<1$. 
    Therefore, the Cartier index of $K_{E_1^\nu}$ at any point on $C_1$ is bounded above by a constant depending only on $I$. 

    \item In this step, we show this Theorem holds when $\varphi_R$, which will be introduced later, is a Mori fibre space. 

    We recall that $K_{E_1^\nu}+\Theta\equiv 0$, $\lfloor\Theta\rfloor=\sum n_iC_i$ with $n_1\geq 2$, $\{\Theta\}=\sum b_i\Theta_i$, the pair $(E_1^\nu,\sum C_i+\sum b_i\Theta_i)$ is log canonical, and no curve $C_i$ with $n_i\geq 2$ intersects $\Theta_1$. 
    
    As $(E_1^\nu,\Theta'-b_1\Theta_1)$ is log canonical, we can run the $(K_{E_1^\nu}+\Theta'-b_1\Theta_1)$-minimal model program. 
    Let $R=\bR_{\geq 0}[G]$ be a $(K_{E_1^\nu}+\Theta'-b_1\Theta_1)$-negative extremal ray where $G$ is an irreducible curve and $\varphi_R \colon E_1^\nu \to S$ be the contraction associated to $R$. 
    We further assume that $G$ intersects $\Theta_1$. 

    Suppose $\varphi_R \colon E_1^\nu \to S$ is a Mori fibre space. 
    Then $S$ is a curve, otherwise the Picard number of $E_1^\nu$ is 1. 
    So $\Theta_1$ is an ample divisor and thus $\Theta_1\cdot C_i>0$ for any $C_i$ with $n_i\geq 2$, which contradicts the Step (a). 
    Thus, $S$ is a curve, and therefore we may assume $G$ is a general fibre of $\varphi_R$. 
    Note that $G$ is nef and 
    \[-K_{E_1^\nu}\cdot G\geq -(K_{E_1^\nu}+\Theta'-b_1\Theta_1)\cdot G>0.\]
    So by Lemma~\ref{fib_neg}, we have $K_{E_1^\nu}\cdot G = -2$ and thus $C_i\cdot G=0$ for all $C_i$ with $n_i\geq 2$. 
    Hence, $b_1\in J_2$ by the choice of $J_2$ and thus $d_1\in I_2$ by the choice of $I_2$. 

    \item In this step, we show the Theorem holds when $\varphi_R$ is a divisorial contraction with exceptional curve $G$ intersecting some $C_i$ with $n_i\geq 2$. 
    
    Note that $G^2<0$ and thus $G\neq \Theta_1$, otherwise, we have 
    \begin{align*}
    0 &> (K_{E_1^\nu}+\Theta'-b_1\Theta_1)\cdot G \\
    &= \Big(-b_1G-\sum (n_i-1)C_i\Big)\cdot G \\
    &= -b_1G^2>0, \mbox{which is impossible.}
    \end{align*}

    Furthermore, we assume $G$ intersects one of $C_i$ with $n_i\geq 2$. 
    For simplicity, let us assume $G$ intersects $C_1$ with $n_1\geq 2$. 
    By adjunction, we have 
    \begin{align*}
    0 > G^2 &= (K_{E_1^\nu}+G+\Theta)\cdot G \\
    &\geq -2 + \frac{m-1+n_1+g}{m}+\frac{n-1+b_1+h}{n} \\
    &\geq -2 + \frac{m-1+2}{m}+\frac{n-1}{n} \\  
    &= \frac{1}{m} - \frac{1}{n} > -1
    \end{align*}
    where $m$ is the Cartier index of $K_{E_1^\nu}$ at some point supported on $C_i\cap G$ and $n$ is the one at some point supported on $\Theta_1\cap G$. 
    Therefore, $n < m\leq M$ by Step (c). 
    Moreover, the inequality above also shows that the Cartier index of $K_{E_1^\nu}$ at any point on $G$ is bounded above by $M$. 

    Since $0>G^2 > -1$ and $(M!)^2G^2$ is an integer, $G^2$ belongs to a finite set depending only on $M$, and therefore depending only on $I$. 
    Thus, by $G^2=(K_{E_1^\nu}+G+\Theta)\cdot G$, adjunction along $G$, and Lemma~\ref{lem_sum_finite}, there is a finite subset $J_5\subset \operatorname{D}(I)$ depending only on $I$ such that $b_1$ belongs to $J_5$. 
    Hence, by \cite[Lemma 5.2]{hacon2014acc}, $d_1$ belongs to a subset $I_5\subset I$, which depends only on $I$. 

    \item In this step, we will demonstrate that we can reduce to the case when all $(K_{E_1^\nu}+\Theta'-b_1\Theta_1)$-negative extremal rays $R=\bR_{\geq 0}[G]$ satisfying $G\cdot \Theta_1=0$. 
    
    So we can assume $\varphi_R\colon E_1^\nu \to S$ is a divisorial contraction such that no curve $C_i$ with $n_i\geq 2$ intersecting $G$. 
    We recall that $G$ intersects $\Theta_1$ by our choice, and thus $G\cdot\Theta_1>0$ as $G\neq\Theta_1$ shown in step (e). 
    
    Note that we have the following properties: 
    \begin{itemize}
        \item $K_S+\Theta_S\equiv 0$ where $\Theta_S$ is the proper transform of $\Theta$ on $S$. 
        \item $(S,\Theta'_S-b_1\Theta_{1,\,S})$ is log canonical as $\Theta_1$ is not contracted by $\varphi_R$ where the subscription $S$ indicates the proper transform of the associated divisors. 
        \item the Cartier index of $K_S$ at any point on $C_{i,\,S}$ with $n_i\geq 2$ is bounded above by $M$ as $G$ does not intersect these $C_i$. 
        \item Besides, $\Theta_{1,\,S}$ does not intersect any $C_{i,\,S}$ with $n_i\geq 2$. 
    \end{itemize} 
    Then we replace $E_1^\nu$ by $S$ and all divisors on $E_1^\nu$ by its proper transform on $S$. 
    We will still use the same notation after replacement. 
    Then we can consider the new contraction $\varphi_{R'}$. 
    By step (d), we use the condition that $\Theta_1$ does not intersect curves $C_i$ with $n_i\geq 2$ to reduce to the case when $\varphi_{R'}$ is a divisorial contraction. 
    By step (e), we use the boundedness of the Cartier index of $K_{E_1^\nu}$ to reduce to the case $\varphi_{R'}$ is a divisorial contraction satisfying no curve $C_i$ with $n_i\geq 2$ intersects $G$. 
    Moreover, all properties listed above are preserved. 
    Therefore, by induction on the Picard number of $E_1^\nu$, we can reduce to the case when all $(K_{E_1^\nu}+\Theta'-b_1\Theta_1)$-negative extremal rays $R=\bR_{\geq 0}[G]$ satisfying $G\cdot \Theta_1=0$. 
    \item In the rest of steps, we demonstrate that the case when all $(K_{E_1^\nu}+\Theta'-b_1\Theta_1)$-negative extremal rays $R=\bR_{\geq 0}[G]$ satisfying $G\cdot \Theta_1=0$ does not happen. 
    In other words, we have already finish the proof of the Theorem. 

    Note that we have following properties:
    \begin{itemize}
        \item $K_{E_1^\nu}+\Theta\equiv 0$. 
        \item $\Theta_1$ does not intersect any curve $C_i$ with $n_i\geq 2$. 
        \item Any $(K_{E_1^\nu}+\Theta'-b_1\Theta_1)$-negative extremal ray $R=\bR_{\geq 0}[C]$ satisfies $C\cdot\Theta_1=0$. 
    \end{itemize}

    As $(K_{E_1^\nu}+\Theta'-b_1\Theta_1)$ is not pseudo-effective and $(E_1^\nu,\Theta'-b_1\Theta_1)$ is log canonical, there is a $(K_{E_1^\nu}+\Theta'-b_1\Theta_1)$-negative extremal ray $R=\bR_{\geq 0}[G]$. 
    By assumption, we have $G\cdot\Theta_1=0$. 

    \item In this step, we show $\varphi_R\colon E_1^\nu\to T$ is not a Mori fibre space. 
    
    Suppose $\varphi_R \colon E_1^\nu \to T$ is a Mori fibre space. 
    Note that $T$ must be a curve, otherwise $\rho(E_1^\nu)=1$ and $G\cdot\Theta_1>0$, which contradicts our assumption.  
    
    Then we may assume $G$ is a general fibre of $\varphi_R$. 
    Since $G\cdot \Theta_1=0$, $\Theta_1$ is in a fibre of $\varphi_R$. 
    We notice that 
    \begin{align*}
    0 &> (K_{E_1^\nu}+\Theta'-b_1\Theta_1)\cdot G \\
    &= \Big(-b_1\Theta_1-\sum (n_i-1)C_i\Big)\cdot G \\
    &= -\sum (n_i-1)C_i\cdot G.
    \end{align*}
    Then one of $C_i$ with $n_i\geq 2$ intersects the general fibre of $\varphi_R$. 
    Let us assume that $C_1$ with $n_1\geq 2$ intersects the general fibre of $\varphi_R$. 
    Thus, $\Theta_1$ is a component of a singular fibre $F$ of $\varphi_R$ with at least two components, otherwise $\Theta_1\equiv aG$ for some $a>0$ and thus $0=\Theta_1\cdot C_1 = aG\cdot C_1 > 0$, which is impossible. 
    Hence, there is another component $C$ of this singular fibre $F$ that intersects $\Theta_1$ because $F$ is connected. 
    Note that 
    \begin{itemize}
    \item $C\neq C_i$ with $n_i\geq 2$ since $C\cdot \Theta_1>0$ but $C_i\cdot\Theta_1=0$ from our assumption, 
    \item $C^2<0$ because $C$ is an irreducible component of a singular fibre of $\varphi_R$, and 
    \item $(K_{E_1^\nu}+\Theta'-b_1\Theta_1)\cdot C = (-b_1\Theta_1-\sum (n_i-1)C_i)\cdot C\leq -b_1\Theta_1\cdot C<0$. 
    \end{itemize}
    So we see that $\bR_{\geq 0}[C]$ is a $(K_{E_1^\nu}+\Theta'-b_1\Theta_1)$-negative extremal ray with $C\cdot\Theta_1>0$, which contradicts our assumption.  
    Thus, we show that $\varphi_R$ is not a Mori fibre space. 

    \item In this step, we show $\varphi_R$ is not a divisorial contraction. 

    We proceed by induction on the Picard number of $E_1^\nu$. 
    It is clear that $\varphi_R$ is not a divisorial contraction when $\rho(E_1^\nu)=1$. 
    Now suppose $\varphi_R\colon E_1^\nu\to T$ is a divisorial contraction. 
    Then we replace $E_1^\nu$ by $T$ and all divisors on $E_1^\nu$ by its proper transform on $T$. 
    We will still use the same notation after replacement. 
    Note that we still have the following properties:
    \begin{itemize}
        \item $K_{T}+\Theta_T\equiv 0$. 
        \item $\Theta_{1,\,T}$ does not intersect any curve $C_{i,\,T}$ with $n_i\geq 2$. 
        \item Any $(K_{T}+\Theta'_T-b_1\Theta_{1,\,T})$-negative extremal ray $R'=\bR_{\geq 0}[G']$ satisfies $G'\cdot\Theta_{1,\,T}=0$. 
    \end{itemize}
    Thus, as $(K_{T}+\Theta'_T-b_1\Theta_{1,\,T}$ is not pseudo-effective, there is a $(K_{T}+\Theta'_T-b_1\Theta_{1,\,T})$-negative extremal ray $R'=\bR_{\geq 0}[G']$ with $G'\cdot\Theta_{1,\,T}=0$. 
    By the proof in step (h), the associated contraction $\varphi_{R'}$ is a divisorial contraction and $\rho(T)=\rho(E_1^\nu)-1$, which contradicts the induction hypothesis.  
    Hence, $\varphi_R$ is not a divisorial contraction. 
\end{enumerate}
\qed
\end{enumerate}
\end{pf}

\section{Foliations of rank one on threefolds}
This section is devoted to show Theorem~\ref{acc_fol_1}, which is Theorem~\ref{thm_coef_finite} when $r=1$ and $n=3$. 

\subsection{Resolution of singularities}
Most definitions in this subsection follow from \cite{cascini2025mmp}. 

\begin{eg}[{\cite[Example III.iii.3]{mcquillan2013almost}}]\label{simple_sing_eg}
    Let $X$ be the quotient of $\bC^3$ by the $\bZ/2\bZ$-action given by $(x,y,z)\mapsto (y,x,-z)$. 
    We consider a vector field on $\bC^3$ given by 
    \[\partial = \Bigg(x\frac{\partial}{\partial x} - y\frac{\partial}{\partial y}\Bigg) + \Bigg(a(xy,z)x\frac{\partial}{\partial x} - a(xy,-z)y\frac{\partial}{\partial y} + c(xy,z)\frac{\partial}{\partial z}\Bigg)\]
    where $a$ and $c$ are formal functions in two variables such that $c$ is not a unit and satisfies $c(xy,z) = c(xy,-z)$. 
    Note that $\partial \mapsto -\partial$ under the $\bZ/2\bZ$-action and thus it induces a foliation $\sF$ on $X$. 
    Moreover, $\sF$ has at worst canonical singularities. 
\end{eg}

\begin{defn}
    Let $X$ be a normal threefold and $\sF$ be rank one foliation with at worst canonical singularities. 
    We say $p\in X$ is a \emph{simple singularity} of $\sF$ if either
    \begin{enumerate}
        \item $\sF$ is terminal and no component of $\operatorname{Sing}(X)$ through $p$ is invariant, 
        \item $X$ and $\sF$ are formally isomorphic to the variety and foliation defined in Example~\ref{simple_sing_eg} at $p$, or
        \item $X$ is smooth at $p$. 
    \end{enumerate}
\end{defn}

\begin{rmk}
    The definition above is for the foliation of \emph{rank one}; while in definition~\ref{defn_simple_corank_one}, we have simple singularity for \emph{co-rank one} foliation. 
    We keep the same terminologies as in the literature. 
    It should be clear which definition we used from the context since one is for co-rank one foliation and another one is for rank one foliation. 
\end{rmk}

\begin{thm}[{\cite[III.iii.4 and III.iii.4.bis]{mcquillan2013almost}}]\label{rsln_n3r1}
    Let $\sF$ be a rank one foliation on a normal threefold $X$. 
    Then there exists a sequence of weighted blowups in foliation invariant centres $\pi \colon \widetilde{X} \to X$ such that $\pi^*\sF$ has only simple singularities. 
    Moreover, there exists a sequence of further blow up along everywhere transverse centres to the induced foliation $\pi'\colon X'\to X$ such that $\pi'^*\sF$ has simple singularities of type (ii) or (iii). 
\end{thm}

\subsection{Adjunction}
\begin{thm}\label{fol_adjunct_rank_1}
    Let $\sF$ be a foliation of rank one with at worst simple singularities on a normal threefold $X$, $S$ be an invariant prime divisor on $X$, and $I$ be a subset of $[0,1]$. 
    Suppose $(\sF,\Delta)$ is log canonical where $\Delta$ is an effective (non-invariant) divisor whose coefficients belong to $I$. 
    
    Let $\nu \colon S^\nu \to S\hookrightarrow X$ be the composition of the normalization and the inclusion. 
    We write $(K_\sF+\Delta)\vert_{S^\nu} = K_\sG+\Theta$ where $\sG$ is the pullback foliation on $S^\nu$ and $\Theta$ is an effective $\bQ$-divisor.  
    Then the coefficients of $\Theta$ belong to $(I\cup\{\frac{1}{2}\})_+\subset(\operatorname{D}(I))_+$. 
\end{thm}
\begin{pf}
    The existence of $\Theta$ and its effectiveness follow from \cite[Proposition 2.13]{cascini2025mmp}. 
    
    Let $C$ be any irreducible component of $\Theta$. 
    If $\sF$ is terminal along $\nu(C)$, then by \cite[Lemma 1.1.3]{bogomolov2016rational}, $\nu(C)$ is non-invariant. 
    Thus, $\nu(C)$ is not contained in $\operatorname{Sing}(X)$ by \cite[Lemma 2.6]{cascini2020mmp}. 
    Hence, $\sF$ is Gorenstein along $\nu(C)$. 
    On the other hand, if $\sF$ is not terminal along $\nu(C)$, then the Cartier index of $K_\sF$ at the generic point of $\nu(C)$ is at most $2$. 
    Moreover, any irreducible component of $\Delta$ does not intersect $S$ along $\nu(C)$ since $\sF$ is not terminal along $\nu(C)$ and $(\sF,\Delta)$ is log canonical. 
    Therefore, the coefficients of $\Theta$ belong to $(I\cup\{\frac{1}{2}\})_+$.
    \qed
\end{pf}

\subsection{ACC for LCT when \texorpdfstring{$n=3$}{n=3} and \texorpdfstring{$r=1$}{r=1}}
\begin{thm}\label{acc_fol_1}
    Fix a set $I\subset [0,1]$, which satisfies the descending chain condition. 
    Then there exists a finite subset $I_0\subset I$ with the following property:
    If $(X,\sF,D)$ is a triple such that 
    \begin{enumerate}
    \item $X$ is a normal threefold,
    \item $\sF$ is a foliation of rank $1$, 
    \item $(\sF,D)$ is log canonical,
    \item the coefficients of $D$ belong to $I$, and 
    \item there is a log canonical centre $Z$ which is contained in every component of $D$,
    \end{enumerate}
    then the coefficients of $D$ belong to $I_0$. 
\end{thm}
\begin{pf}
Without loss of generality, we may assume that $Z$ is maximal with respect to the inclusion. 
If $Z$ is a divisor, then $Z = D$ and thus we require $1\in I_0$. 

Now we suppose $\dim Z\leq 1$. 
By the resolution of singularities (Theorem~\ref{rsln_n3r1}) along $Z$, we have a birational morphism $p \colon X' \to X$ such that the pullback foliation $\sF':= p^*\sF$ has only simple singularities. 
Let $D' := p_*^{-1}D$. 
Note that $\sF$ is terminal along $Z$ by \cite[Lemma 8.3]{cascini2025mmp}. 
Thus, all exceptional divisors with centre $Z$ are invariant. 
Then we may write 
\[K_{\sF'}+D'= p^*(K_\sF+D)+F'\]
where $F'$ is an effective $p$-exceptional invariant divisor. 
After possibly taking a higher resolution, we may assume that $(\sF',D')$ is log canonical. 
By \cite[Proposition 8.1]{cascini2025mmp}, we may run a $(K_{\sF'}+D')$-MMP over $X$. 
Let $\varphi \colon X' \dashrightarrow Y$ be the output of this MMP with $\pi \colon Y \to X$. 
Let $\sG := \varphi_*\sF'$ and $D_Y := \varphi_*D'$. 
Then we have, by \cite[Proposition 8.1 and Lemma 8.3]{cascini2025mmp}, that 
\begin{enumerate}
\item $(\sG,D_Y)$ is log canonical;
\item we have $K_\sG + D_Y = \pi^*(K_{\sF} + D)$;
\item $\sG$ has only simple singularities; and
\item $\varphi$ contracts exactly those exceptional divisors with positive discrepancies. 
\end{enumerate}

Let $D_1$ be one of irreducible components of $D$ and $d_1$ be the coefficient of $D_1$ in $D$. 
Then the proper transform $\pi_*^{-1}D_1$ has coefficient $d_1$ in $\Gamma$. 
Note that there is at least one irreducible component of exceptional divisor $E_j$ intersecting $\pi_*^{-1}D_1$. 
Without loss of generality, we may assume that $E_1$ intersects $\pi_*^{-1}D_1$. 
Let $\Theta_1$ be an irreducible component of $\pi_*^{-1}D_1\cap E_1$. 
By the assumption that $Z$ is contained in $D_1$, we have $\Theta_1$ dominates $Z$. 

By Theorem~\ref{fol_adjunct_rank_1}, we have 
\[(K_\sG+D_Y)\vert_{E_1^\nu} = K_\sH+\Theta\]
where $\sH$ is the pullback foliation on $E_1^\nu$ and $\Theta$ is an effective divisor whose coefficients belong to $(I\cup\{\frac{1}{2}\})_+$. 
Let $\Theta = \sum b_i\Theta_i$ where $\Theta_1$ be an irreducible component of $\pi_*^{-1}D_1\cap E_1$ which dominates $Z$. 
Note that $b_1 = kd_1+f$ where $k\in\bN$ and $f\in (I\cup\{\frac{1}{2}\})_+\subset(\operatorname{D}(I))_+$. 
\begin{enumerate}
\item If $Z$ is a curve, then $\psi := \pi\vert_{E_1^\nu}\colon E_1^\nu \to Z$ is a fibration. 

By the assumption that $Z$ is contained in $D_1$, we have $\Theta_1\cdot F > 0$ where $F$ is a general fibre of $\psi$. 
Note that  
\begin{equation}\label{r1ve0curve}
0 = \pi^*(K_\sF+D)\cdot F = (K_\sH+\Theta)\cdot F
\end{equation}
and thus $-K_\sH\cdot F = \Theta\cdot F\geq \Theta_1\cdot F >0$. 
By Lemma~\ref{fib_neg}, we have $F$ is invariant and $K_\sH\cdot F = -2$. 
Therefore, we have
\[2 = -K_\sH\cdot F = \Theta\cdot F = \sum b_i(\Theta_i\cdot F)\]
where $\Theta_i\cdot F \in\bZ_{\geq 0}$ and $\Theta_1\cdot F>0$. 
Replacing $b_1$ by $kd_1+f$, we get 
\[2 = (\Theta_1\cdot F)kd_1 + g\]
where $g\in(I\cup\{\frac{1}{2}\})_+$. 
Since $I\cup\{\frac{1}{2}\}$ satisfies the descending chain condition, by Lemma~\ref{lem_sum_finite}, there is a finite subset $J_1\subset I\cup\{\frac{1}{2}\}$ such that $d_1\in J_1$. 
Thus, $d_1$ belongs to $I_1 := J_1\cap I$ which a finite subset of $I$. 

\item Now suppose $Z$ is a point. 
Then $K_\sH+\Theta\equiv 0$. 
We will divide the proof into the following steps:

\begin{enumerate}
    \item In this step, we show $\sH$ is algebraically integrable. 

    Let $A$ be an ample divisor on $E_1^\nu$. 
    Since $\Theta$ is effective and $\Theta_1\subset\operatorname{Supp}(\Theta)$, we have $\Theta\cdot A>0$ and thus, $K_\sH\cdot A<0$. 
    Hence, the slope of $\sH$ with respect to $A$ is $\mu_A(\sH) = -K_\sH\cdot A > 0$. 
    Since $\sH$ is of rank one, we have $\mu_{A,\,\textnormal{min}}(\sH) = \mu_A(\sH)>0$. 
    By \cite[Theorem 1.1]{campana2019foliations} or \cite[Proposition 2.2]{ou2023generic}, $\sH$ is algebraically integrable. 

    \item Let $L$ be a general leaf of $\sH$ and $\nu\colon E_1^\nu\to E_1\hookrightarrow Y$ be the composition of the normalization and the inclusion. 
    Note that $\sG$ is terminal along $\nu(\Theta_1)$ as $\nu(\Theta_1)\subset\pi_*^{-1}D_1$, and thus $\nu(\Theta_1)$ is not invariant by \cite[Lemma 1.1.3]{bogomolov2016rational}. 
    Therefore, $\Theta\cdot L>0$ and $D_Y\cdot\nu(L)>0$. 

    \item Suppose there is no dicritical singularity for $\sH$. 
    As $L$ is general, we may assume $L$ is smooth and disjoint from the singular locus of $E_1^\nu$. 
    Thus, we have 
    \[0 = (K_\sH+\Theta)\cdot L = -2+\Theta\cdot L.\] 
    Therefore, 
    \[2 = -K_\sH\cdot L = \Theta\cdot L = \sum b_i(\Theta_i\cdot L)\]
    where $\Theta_i\cdot L\in\bZ_{\geq 0}$ and $\Theta_1\cdot L>0$. 
    Replacing $b_1$ by $kd_1+f$, we get 
    \[2=(\Theta_1\cdot L)kd_1+g\]
    where $g\in (I\cup\{\frac{1}{2}\})_+$. 
    Then by the construction of $I_1$, we have $d_1\in I_1$. 

    \item Now we suppose there is a dicritical singularity $p$ for $\sH$. 
    In this step, we show that $\nu(p)$ is a log canonical centre for $\sG$, as well as a singularity for $\sG$.  

    Assume $\nu(p)$ is not a log canonical centre for $\sG$.  
    Since $\sG$ has simple singularities which are log canonical, $\sG$ is terminal at $\nu(p)$. 
    Note that $K_\sG\cdot \nu(L) = -D_Y\cdot \nu(L) <0$ for a general leaf $L$ of $\sH$. 
    By \cite[Proposition 3.3 and Theorem 3.2]{cascini2025mmp}, $\sG$ is a fibration locally around $\nu(p)$, which contradicts that $p$ is a dicritical singularity for $\sH$. 
    Therefore, $\nu(p)$ is a log canonical centre for $\sG$. 
    Moreover, $\nu(p)$ is a singularity for $\sG$ by \cite[Lemma 2.9]{cascini2025mmp}. 

    \item In this step, we will finish the proof for the case when $\sH$ has dicritical singularities. 
    
    By \cite[Proposition 2.13]{cascini2025mmp}, we have $K_\sG\vert_{L^\nu}=K_{L^\nu}+\Delta$ where the image of $\operatorname{Supp}\lfloor\Delta\rfloor$ on $Y$ contains all singularities of $\sG$ on $\nu(L)$. 
    Since $0>K_\sG\cdot \nu(L) = -2+\deg\Delta$, there is at most one singularity for $\sG$. 
    As $\sH$ has dicritical singularities, whose images on $Y$ are singularities for $\sG$, there is exactly one dicritical singularity, say $p$, and there is no singularity for $\sG$ on $\nu(L)$ except for $\nu(p)$. 
    
    As $\sG$ has only simple singularities, $2K_\sG$ is Cartier at $\nu(p)$ and thus, the coefficient of $\Delta$ at $p$ is $1$ or $\frac{3}{2}$, both of which belong to $(\operatorname{D}(\{0\}))_+\subset(\operatorname{D}(I))_+$. 
    Moreover, by \cite[Proposition 2.13]{cascini2025mmp}, the coefficients of $\Delta$ at points other than $p$ belong to $\operatorname{D}(\{0\})\subset(\operatorname{D}(I))_+$. 

    On the other hand, we note that $D_Y\vert_{E_1^\nu}\in I_+$ by Theorem~\ref{fol_adjunct_rank_1}. 
    As $\nu(p)$ is a log canonical centre for $\sG$, $\operatorname{Supp}(D_Y)$ does not contain $\nu(p)$. 
    Moreover, as $L$ is general, $L$ is smooth outside $p$ and thus, $D_Y\cdot \nu(L) = D_Y\vert_{E_1^\nu}\cdot L = \ell b_1+g$ where $\ell\in\bN$ and $g\in I_+\subset(\operatorname{D}(I))_+$. 

    To sum up, from $0=(K_\sH+\Theta)\cdot L = (K_\sG+D_Y)\cdot \nu(L)$, we get 
    \[2 = \deg\Delta+\ell b_1+g = \deg\Delta+\ell(kd_1+f)+g\]
    where both $\deg\Delta$ and $f$ belong to $(\operatorname{D}(I))_+$. 

    Note that $\operatorname{D}(I)$ satisfies the descending chain condition. 
    By Lemma~\ref{lem_sum_finite}, there is a finite subset $J_2\subset\operatorname{D}(I)$ depending only on $I$ such that $d_1\in I_2$. 
    Thus, $d_1$ belongs to $I_2 := J_2\cap I \subset I$, which is a finite subset of $I$. 
\end{enumerate}
\qed
\end{enumerate}
\end{pf}

\section{Proof of Theorem on ACC for LCT}
\begin{proof}[{Proof of Theorem~\ref{main_acc}}]
    Let $\{c^{(\ell)}\}_{\ell=1}^\infty$ be an increasing sequence in $\operatorname{LCT}_{n,r}(I,J)$. 
    We may assume that both 0 and 1 belong to $I$ and $J$. 
    
    For each $c^{(\ell)}$ in the sequence, we may find a foliation $\sF^{(\ell)}$ of rank $r$ on a variety $X^{(\ell)}$ of dimension $n$, a divisor $\Delta^{(\ell)}$ whose coefficients belong to $I$, and an $\bR$-Cartier divisor $M^{(\ell)}$ whose coefficients belong to $J$ such that $c^{(\ell)} = \sup\{t\in\bR\vert\, (\sF^{(\ell)},\Delta^{(\ell)}+tM^{(\ell)}) \textnormal{ is log canonical}\}$. 
    
    Let $D^{(\ell)} = \Delta^{(\ell)}+c^{(\ell)}M^{(\ell)}$ and $K = \{a+bc^{(\ell)}\leq 1\vert\,a\in I, b\in J, \textnormal{ and for any $\ell$}\}$. 
    Note that $K$ satisfies the descending chain condition and $(\sF^{(\ell)},D^{(\ell)})$ is log canonical with coefficients of $D^{(\ell)}$ in $K$. 
    Since $c^{(\ell)}$ is the log canonical threshold, there is a log canonical centre $Z^{(\ell)}$ contained in the support of $M^{(\ell)}$. 
    Possibly throwing away components of $D^{(\ell)}$ which do not contain $Z^{(\ell)}$, then we may assume every component of $D^{(\ell)}$ contains $Z^{(\ell)}$. 
    By Theorem~\ref{thm_coef_finite}, there is a finite subset $K_0\subset K$ such that the coefficients of $D^{(\ell)}$ belong to $K_0$. 
    This shows that the sequence $c^{(\ell)}$ eventually terminates. 
\end{proof}

\bibliographystyle{amsalpha}
\bibliography{ACC_LCT_Final}

\end{document}